\documentclass{article}
\usepackage{amsthm}
\usepackage[sumlimits]{amsmath}
\usepackage{amsfonts}
\usepackage{color, framed, url}
\usepackage{graphicx}
\usepackage[margin=2cm]{geometry}
\DeclareMathOperator{\ad}{ad}

\DeclareMathOperator{\diff}{d\!}
\DeclareMathOperator{\Emb}{Emb}
\DeclareMathOperator{\Diff}{Diff}

\newtheorem{theorem}{Theorem}
\newtheorem{lemma}[theorem]{Lemma}
\newtheorem{corollary}[theorem]{Corollary}
\newtheorem{definition}[theorem]{Definition}

\newtheorem{remark}[theorem]{Remark}
\newcommand{\rem}[1]{}
\def\MM#1{\boldsymbol{#1}}
\newcommand{\pp}[2]{\frac{\partial #1}{\partial #2}} 
\newcommand{\dede}[2]{\frac{\delta #1}{\delta #2}}
\newcommand{\iprod}[3]{\left\langle #1,#2 \right\rangle_{#3}}
\newcommand{\dd}[2]{\frac{\diff#1}{\diff #2}}

\newcommand{\bfi}[1]{{\bfseries\itshape #1}}

\begin{document}

\title{Geodesic boundary value problems with symmetry} 
\author{C. J. Cotter$^{1}$ and D. D. Holm$^{2}$}
\addtocounter{footnote}{1}
\footnotetext{Department of Aeronautics, Imperial College London. London SW7 2AZ, UK. 
\texttt{colin.cotter@imperial.ac.uk}
\addtocounter{footnote}{1} }
\footnotetext{Department of Mathematics, Imperial College London. London SW7 2AZ, UK. Partially supported by Royal Society of London Wolfson Award.
\texttt{d.holm@imperial.ac.uk}
\addtocounter{footnote}{1} }

\date{December 21, 2009}
\maketitle

\makeatother

\maketitle



\begin{abstract}
  This paper shows how commuting left and right actions of Lie groups
  on a manifold may be used to complement one another in a variational
  reformulation of optimal control problems as geodesic boundary value
  problems with symmetry. In such problems, the endpoint boundary
  condition is only specified up to the right action of a symmetry
  group. In this paper we show how to reformulate the problem by
  introducing extra degrees of freedom so that the endpoint condition
  specifies a single point on the manifold. We prove an equivalence
  theorem to this effect and illustrate it with several examples.  In
  finite-dimensions, we discuss geodesic flows on the Lie groups
  $SO(3)$ and $SE(3)$ under the left and right actions of their
  respective Lie algebras.  In an infinite-dimensional example, we
  discuss optimal large-deformation matching of one closed curve to
  another embedded in the same plane. In the curve-matching example,
  the manifold $\Emb(S^1, \mathbb{R}^2)$ comprises the space of closed
  curves $S^1$ embedded in the plane $\mathbb{R}^2$.  The
  diffeomorphic left action $\Diff(\mathbb{R}^2)$ deforms the curve by
  a smooth invertible time-dependent transformation of the coordinate
  system in which it is embedded, while leaving the parameterisation
  of the curve invariant. The diffeomorphic right action $\Diff(S^1)$
  corresponds to a smooth invertible reparameterisation of the $S^1$
  domain coordinates of the curve. As we show, this right action
  unlocks an important degree of freedom for geodesically matching the
  curve shapes using an equivalent fixed boundary value problem,
  without being constrained to match corresponding points along the
  template and target curves at the endpoint in time.

\end{abstract}

\section{Introduction}
In this paper we are concerned with finding geodesics between points
on manifolds. The construction of geodesics is useful for studying
problems on manifolds since they can describe the relationship between
two points. Within a coordinate patch on a manifold, any point can
described relative to a reference point by specifying a direction and
a length along the geodesic in that direction. This becomes useful for
performing statistics on the coordinate patch, for example.  In this
paper we consider problems in which the endpoint of the trajectory is
only fixed up to the orbit of a Lie group. In low dimensional cases
(and we shall describe some examples of these) it is often easy to
solve these problems by constructing reduced coordinates which do not
change under the action of the Lie group. However, in many cases it is
difficult to construct such coordinates, especially if the problem is
to be discretised and solved numerically. In this paper we provide a
framework that allows one to work with full unreduced coordinates on
the manifold, by transforming to an equivalent problem which has the
endpoint of the trajectory fixed exactly.\\[1mm]

There are many examples of problems where this framework can be
applied, but we are motivated by the problem of obtaining
diffeomorphisms on $\mathbb{R}^2$ which map one embedded curve
$\Gamma^A$ into another embedded curve $\Gamma^B$, and which minimise
a given metric so that they are geodesics in the diffeomorphism
group. The aim is to find a characterisation of curve $\Gamma^B$ with
respect to curve $\Gamma^A$ that is independent of parameterisations
of the curves. This means that we do not specify \emph{a priori} the
point on $\Gamma^A$ which gets matched to each specific point on
$\Gamma^B$, and so the minimisation is performed over all
parameterisations of the curves. In practise the computation is
performed using a particular parameterised curve $\MM{q}\in
\Emb(S,\mathbb{R}^2)$ (where $S$ is the embedded space, for example,
the circle for simple closed curves). In computing the equations of
motion, a conjugate momentum $\MM{p}_q\in T^*_{q}\Emb(S,\mathbb{R}^2)$
is constructed, and the flow taking the initial curve $\Gamma^A$ to
the final curve $\Gamma^B$ can be characterised entirely by the
initial conditions $\MM{p}_q|_{t=0}$ for the conjugate momentum. In
fact, it turns out that $\MM{p}_q|_{t=0}$ is normal to the curve, so the
flow can be characterised by a one-dimensional signal. Since
$T^*_q\Emb(S,\mathbb{R}^2)$ is a linear space, linear statistics can
be computed on $\MM{p}_q|_{t=0}$. For example, this may allow one to
test the hypothesis that there is a statistical correlation between
between the shape of the surface of a biological organ, obtained from
a medical scan, and future development of disease.\\

To discuss the issues further, we formulate the curve matching problem
described above, which may be regarded as an optimal control problem
in the sense of the problems discussed in
\cite{BlCrHoMa2000,BlCrMaRa1998}:
\begin{definition}[\bfi{Curve matching problem}]
\label{curve matching problem}
  Let $\MM{q}(s;t)$ be a one-parameter family of parameterised simple closed 
  curves in $\mathbb{R}^2$, with $s\in [0,1]$ being the curve parameter and
  $t \in [0,1]$ being the parameter for the family. Let
  $\MM{u}(\MM{x};t)$ be a one-parameter family of vector fields on
  $\mathbb{R}^2$. Let $\eta$ be a diffeomorphism of $S^1$. We seek
  $\MM{q}$ and $\MM{u}$ which satisfy
\[
\min_{\MM{u},\eta}\int_0^1\frac{1}{2}\|\MM{u}\|^2_{V}\diff{t}
\]
subject to the constraints
\begin{eqnarray} \mbox{\rm [Reconstruction relation]} \quad 
\pp{}{t}\MM{q}(s;t)
  &=& \MM{u}(\MM{q}(s;t),t), \label{dynamical constraint1}  \\
\mbox{\rm Initial state (Template)]} \quad 
\MM{q}(s;0) & = & \MM{q}^A(s), \label{template constraint1} \\
\mbox{\rm [Final state (Target)]} \quad 
\MM{q}(s;1) & = & \MM{q}^B(\eta(s)), \label{target constraint1}
\end{eqnarray}
where $\|\cdot\|_{V}$ is the chosen norm which defines the space of
vector fields $V$.
\end{definition}
The solution of this problem describes a geodesic in the
diffeomorphism group which takes the simple closed curve $\Gamma^A$
parameterised by $\MM{q}^A$ to the simple closed curve $\Gamma^B$
parameterised by $\MM{q}^B$. We represent the shapes of simple closed
curves as elements of $\Emb(S^1,\mathbb{R}^2)/\Diff(S^1)$, where
$\Diff(S^1)$ is the group of diffeomorphisms of $S^1$. However, 
we do not want to calculate on this space; instead, we want to calculate
on the full space $\Emb(S^1,\mathbb{R}^2)$ by minimising over
all reparameterisations $\eta(s)\in\Diff(S^1)$. \medskip

There are two general strategies for solving such problems. The first
strategy, used for example in \cite{CaYo01}, is to use a gradient
method (\emph{i.e.} a modification of the steepest descent method such
as the nonlinear conjugate gradient method \cite[and references
therein]{Sh1994}) to minimise the action integral over paths
$\MM{q}(s,t)$ which satisfy the dynamical constraint (this constraint
was enforced ``softly'' \emph{via} a penalty term in
\cite{CaYo01}). An alternative method, referred to in
\cite{MiMaSh2006} as the ``Hamiltonian method'', is to introduce
Lagrange multipliers $\MM{p}(s,t)$ which enforce the dynamical
constraint, and to derive Hamilton's canonical equations for $\MM{q}$ and
$\MM{p}_q$, following the general derivation described in \cite[for
example]{CoHo2009}. Minimisation over the reparameterisation $\eta$,
together with a conservation law obtained from Noether's theorem,
results in the condition that the tangential component of $\MM{p}_q$
vanishes. The aim of the Hamiltonian method is to turn an optimisation
problem into an algebraic equation given by the time-1 flow map of Hamilton's canonical equations. One then solves a shooting problem to find
initial conditions for the normal component $\MM{p}_q$ which generate
solutions to Hamilton's equations that satisfy the boundary
condition \eqref{target constraint1}. The difficulty in solving this
problem numerically lies in finding a good numerical discretisation of
the target constraint condition \eqref{target constraint1}. Various
functionals have been proposed which vanish when the constraint
condition is satisfied. In \cite{GlTrYo04} a functional was proposed
based on singular densities (measures), and in \cite{VaGl2005} a
functional was proposed based on singular vector fields (currents). An
alternative spatial discretisation for the current functional based on
particle-mesh methods was proposed in \cite{Co2008}. There are several
difficulties with these functionals: one is that after numerical
discretisation the functionals do not vanish at the minima, and the
boundary condition must be replaced by a functional minimising
condition. It is also difficult to express the probability
distribution of the functional given the distribution of measurement
errors; this is important for statistical modelling.  \medskip

In
this paper we consider a transformation of problems of the above type,
which results in an alternative formulation that removes the
reparameterisation variable $\eta$ from the target constraint,
thereby resulting in a standard two point boundary value problem on
$T^*\Emb(S,\mathbb{R}^2)$ (with a constraint on the initial
conditions plus an additional parameter). This transformation can be
applied to a very general class of problems; so we present it in
the general case of Lie group actions on a manifold. \medskip

The rest of this paper is organised as follows. In Section
\ref{general section}, we formulate the optimal control problem, then transform to the geodesic problem with symmetry and prove that the two problems are equivalent. In Section \ref{examples} we give some examples and discuss the application to matching curves and surfaces. Section \ref{summary} is the summary and outlook.

\section{Reparameterised geodesic boundary value problems with
  symmetry}
\label{general section}
In this section we describe a general framework for geodesic boundary
value problems with symmetry. We define the following \bfi{Optimal
  Control Problem}.
\begin{definition}[\bfi{Geodesic boundary value problem with symmetry}]
\label{matching problem}
Let $Q$ be a manifold, let $G$ be a Lie group acting on $Q$ from the
left, and let $H$ be a (possibly different) Lie group acting on $Q$
from the right that commutes with the left action of $G$ on $Q$, with
corresponding Lie algebras $\mathfrak{g}$ and $\mathfrak{h}$, and
corresponding Lie algebra actions $\mathcal{X}^G$ and $\mathcal{X}^H$
respectively.  Furthermore, let
$\mathcal{A}:\mathfrak{g}\to\mathfrak{g}^*$ be a positive-definite
self-adjoint operator and let $\langle
\cdot,\cdot\rangle_{\mathfrak{g}}:
\mathfrak{g}\times\mathfrak{g}^*\to\mathbb{R}$ be a nondegenerate
pairing which defines an inner product on $\mathfrak{g}$. We seek
\begin{itemize}
\item a one parameter family ${q}$ of points on $Q$ parameterised by
 $t\in[0,1]$,
\item a one parameter family ${\xi}$ of elements of
  $\mathfrak{g}$ for $t\in[0,1]$, and
\item $\eta\in H$,
\end{itemize}
which minimise
\[
\int_0^1\frac{1}{2}\langle{\xi},\mathcal{A}{\xi}\rangle
_{\mathfrak{g}}\diff{t}
\,,
\]
subject to the constraints
\begin{eqnarray} \mbox{\rm [Reconstruction relation]} \quad 
\dd{}{t}{q}
  &=& \mathcal{X}^G_{{\xi}}{q}, \label{dynamical constraint}  \\
\mbox{\rm [Initial state (Template)]} \quad 
{q}|_{t=0} & = & {q}^A, \label{template constraint} \\
\mbox{\rm [Final state (Target)]} \quad 
{q}|_{t=1} & = & R_{\eta}{q}^B, \label{target constraint}
\end{eqnarray}
where ${q}^A$, ${q}^B$ are chosen points on $Q$, and $R_\eta$ is the
right-action of $\eta$ on $Q$.
\end{definition}
\begin{remark}
  This problem is an optimal control problem in which we seek the
  shortest path in $Q$ from ${q}^A$ to any point ${q}^B\eta$, $\eta\in
  H$. This means we are seeking the shortest path in $Q/H$, but are
  performing the computation on $Q$. In many cases it is much easier
  to compute on $Q$, for example when $Q$ is a vector space. We refer
  to this process of solving a problem on $Q/H$ by calculating on $Q$
  as ``un-reduction''.
\end{remark}

One approach to solving this problem is to derive equations of motion
for $q$, $\xi$ and an optimal condition for $\eta$ and then solving a
shooting problem to find $\eta$ and the initial conditions for $\xi$
which allow equation \eqref{target constraint} to be satisfied.  We
can derive the equations of motion by enforcing the reconstruction
relation \eqref{dynamical constraint} as a constraint using Lagrange
multipliers $p_q\in T_{{q}}^*Q$. This approach leads to the following
variational principle.

\begin{definition}[\bfi{Variational principle for geodesic boundary
    value problem with symmetry}]
\label{variational principle matching problem}
We seek $({p},q)\in T^*Q$ and $\xi\in\mathfrak{g}$ for
$t\in[0,1]$, and $\eta\in H$, which satisfy
\begin{eqnarray}
\delta S = \delta \int_0^1
\frac{1}{2}\left\langle \xi, \mathcal{A}\xi
\right\rangle_{\mathfrak{h}}
+ \left\langle p_q,\dd{}{t}q - \mathcal{X}^G_{\xi}q \right\rangle
_{T^*Q}
\diff{t} = 0,
\end{eqnarray}
subject to 
\begin{eqnarray}
q|_{t=0} = q^A, \quad q|_{t=1} = R_{\eta}q^B,
\end{eqnarray}
where we allow $p_q$, $q$, $\xi$ and $\eta$ to vary.
\end{definition}
From this variational principle we can derive the equations of motion,
which can be used in solving the shooting problem. Before we do this,
we recall the definition of the cotangent-lifted momentum map:
\begin{definition}
  Given an action of a Lie algebra $\mathfrak{g}$ on $Q$, the
  cotangent-lifted momentum map $\mathbf{J}:T^*Q\to \mathfrak{g}$ is defined
  from the formula
\begin{eqnarray}
  \left\langle  \mathbf{J}(p_q),\, \zeta \right\rangle
  _{\mathfrak{g}}
  = 
  \left\langle p_q, \mathcal{X}_{\zeta}q \right\rangle
  _{T^*Q}
\label{diamond-def}
\end{eqnarray}
for all $\zeta\in\mathfrak{g}$. Since we have two Lie algebra actions,
we shall write $\mathbf{J}_G$ for the cotangent-lifted momentum map
corresponding to the left action $\mathcal{X}^G$ of $\mathfrak{g}$ on
$Q$, and $\mathbf{J}_H$ for the cotangent-lifted momentum map corresponding to
the right action $\mathcal{X}^H$ of $\mathfrak{h}$ on $Q$.
\end{definition}

\begin{lemma}[\bfi{Equations of motion for geodesic problem}]
\label{matching problem eqns}
At the optimum, the following equations are satisfied (weakly, for appropriate pairings):
\begin{eqnarray}
\label{dot q}
\dd{}{t}q - \mathcal{X}^G_\xi q & = & 0, \\
\label{dot p}
\dd{}{t}p_q + 
\left(T_{q}\left(\mathcal{X}_\xi^G q\right)\right)^*p_q & = & 0, \\
\label{xi}
\mathcal{A}\xi - \mathbf{J}_G(p_q) & = & 0.
\end{eqnarray}
Furthermore,
\begin{equation}
\label{end condition}
\mathbf{J}_H(p_q)|_{t=1} = 0.
\end{equation}
\end{lemma}
\begin{remark}
  The end-point condition \eqref{end condition} at time t=1 arises from
  minimising over $\eta$ and ensures we have the shortest path over
  $Q/H$.
\end{remark}
\begin{proof} The proof is a direct calculation. 
\begin{eqnarray*}
\delta S 
& = & 
\int_0^1 \iprod{\delta \xi}{\mathcal{A}\xi}{\mathfrak{g}}
+\iprod{\delta p_q}{\dd{q}{t}}{T^*Q}
+\iprod{p_q}{\delta\dd{q}{t}-\delta(\mathcal{X}_{\xi}q)}{T^*Q}\diff{t} \\
 & = & \int_0^1 \iprod{\delta \xi}{\mathcal{A}\xi - \mathbf{J}_G(p_q)}{\mathfrak{g}}
+ \iprod{\delta p_q}{\dd{q}{t}-\mathcal{X}_{\xi}q}{T^*Q}
- \iprod{\dd{p_q}{t} + \left(T_{q}\left(\mathcal{X}_\xi^G q\right)\right)^*p_q}
{\delta q}{T^*Q} \diff{t} + \left[
\iprod{p_q}{\delta q}{T^*Q}
\right]_{t=0}^1.
\end{eqnarray*}
Since $\delta p$, $\delta q$ and $\delta \xi$ are all arbitrary,
stationarity $\delta S=0$ implies equations (\ref{dot q}-\ref{xi}) and their appropriate pairings. The boundary term becomes
\begin{eqnarray*}
\left[
\iprod{p}{\delta q}{T^*Q}
\right]_{t=0}^1 & = & \iprod{p}
{T_{\eta}\left(R_{\eta}q\right)\cdot \delta \eta}{T^*Q}\Big|_{t=1} \\
&=&\iprod{p}{\mathcal{X}^H_{\gamma}q}{T^*Q}\Big|_{t=1} \\
&=&\iprod{\mathbf{J}_H(p_q)}{\gamma}{\mathfrak{h}}\Big|_{t=1} ,
\end{eqnarray*}
where $\gamma$ is the generator of $\delta \eta$, 
\emph{i.e.}
\[
\delta \eta = \dd{}{\epsilon}\Big|_{\epsilon=0}\exp(\epsilon\gamma)\eta.
\]
Hence, we also obtain equation \eqref{end condition} and its appropriate pairing.
\end{proof}
Lemma \ref{matching problem eqns} reformulates the geodesic calculation as a shooting problem in which
one seeks initial conditions for $p_q$ such that $q|_{t=1}=R_\eta q^B$
where $\eta$ is fixed by the condition \eqref{end condition}. Next we
show conservation of the momentum map $\mathbf{J}_H$; this will enable us to
transfer condition \eqref{end condition} from $t=1$ to $t=0$.
\begin{lemma}[\bfi{Noether's theorem for geodesic problem}]
\label{noether matching}
The system of equations (\ref{dot p}-\ref{xi}) has a conserved
momentum $\mathbf{J}_H(p_q)$.
\end{lemma}
\begin{proof}
  The problem in Definition \ref{variational principle matching
    problem} is invariant under transformations
\[
q \to R_{\alpha}q, \qquad \alpha\in H,
\]
which are generated by $\gamma\in\mathfrak{h}$. This means that the
variational principle in Definition \ref{variational principle
  matching problem} is invariant under application of the
\bfi{cotangent lift} (\emph{i.e.}, the dual of the inverse of its
infinitesimal transformation in $Q$), namely
\[
\delta q = \mathcal{X}^H_{\gamma}q, \quad \delta p_q =
-\left(T_{q}\left( \mathcal{X}^H_{\gamma}q\right)\right)^*p_q,
\quad \delta\omega =0,
\]
where for convenience we have inserted the time dependence
\[
\gamma = 0 \quad\mbox{if}\quad t_0<t<1
.\]
 Substitution of this
infinitesimal transformation into the variational principle gives
\begin{eqnarray*}
0 = \delta S 
& = & 
\int_0^{t_0} 
\left[ \iprod{\delta p_q}{\dd{}{t}q - \mathcal{X}^G_{\xi}q}{T^*Q}
+ 
\iprod{p_q}{\delta \dd{}{t}q - T_{q}\left(\mathcal{X}^G_{\xi}q\right)
\cdot\delta q}{T^*Q}
\right]\diff{t}
\\
& = & \int_0^{t_0}
\left[
\iprod{\delta p_q}{\underbrace{\dd{}{t}q - \mathcal{X}^G_{\xi}q}_{=0}}{T^*Q}
-\iprod{\underbrace{\dd{}{t}p_q + \left(T_{q}\left(\mathcal{X}^G_{\xi}\right)
\right)^*p_q}_{=0}}{\delta q}{T^*Q}
\right]\diff{t} + \left[
\iprod{p_q}{\delta q}{T^*Q}
\right]_{t=0}^{t=t_0} \\
&=& \left[
\iprod{p_q}{\mathcal{X}^H_{\gamma}q}{T^*Q}
\right]_{t=0}^{t=t_0} \\
&=& \left[
\iprod{\gamma}{\mathbf{J}_H(p_q)}{\mathfrak{h}}\right]_{t=0}^{t=t_0}. \\
\end{eqnarray*}
Since this equation holds for any $0<t_0<1$, it follows that the quantity 
$\mathbf{J}_H(p_q)$ is conserved.
\end{proof}
Combining this conservation result with equation \eqref{end condition} gives the following
easy corollary.
\begin{corollary}[\bfi{Vanishing  momentum}]
\label{vanishing momentum}
The conserved  momentum satisfies $J(p_q)=0$ for all times $0\leq t
\leq 1$.
\end{corollary}
\begin{proof}
  Lemma \ref{matching problem eqns} states that this quantity vanishes
  for $t=1$, and Lemma \ref{noether matching} states that it is conserved; hence, it always vanishes.
\end{proof}
Corollary \ref {vanishing momentum} implies that solutions of the optimal control problem all have vanishing right action momentum map $\mathbf{J}_H(p_q)=0$. This is what facilitates the ``un-reduction''. Namely, we can compute on $Q$ instead of
$Q/H$ by keeping $\mathbf{J}_H(p_q)=0$.  To obtain the shortest path between
two points in $Q/H$ by solving in $Q$, select a point $q\in Q$ which
is a member of the equivalence class which is the initial point in
$Q/H$, and find initial conditions for $p_q$ such that $\mathbf{J}_H(p_q)=0$; so
that the solution to equations (\ref{dot q}-\ref{xi}) satisfies
$q|_{t=1}=R_\eta q$ for some $\eta\in H$. Computationally, there are
reasons why solving the problem in this form may be difficult. In Section \ref{sec: curve matching}, we
shall describe how the difficulty arises for the curve matching problem specified in the Introduction. In this paper, we shall introduce a reformulation of
the problem for which there is a single fixed value of $q|_{t=1}$.

Before introducing the reformulation, we define the $\operatorname{ad}$ and $\operatorname{ad}^*$ operations for the  Lie algebra $\mathfrak{g}$ and briefly discuss the reduced equation for the Lie algebra variable $\xi\in\mathfrak{g}$. The latter is the Euler-Poincar\'e equation for Hamilton's principle with Lagrangian given by the 
energy $\langle \xi,\mathcal{A}\xi \rangle_{\mathfrak{g}}/2$, where $\mathcal{A}:\mathfrak{g}\to\mathfrak{g}^*$ is the positive-definite
self-adjoint operator in Definition \ref{matching problem} of the geodesic matching problem.

\begin{definition}[\bfi{Notation for the $\operatorname{ad}$ and $\operatorname{ad}^*$ operations}]
We define the operation $\operatorname{ad}:\mathfrak{g}\times\mathfrak{g}\to\mathfrak{g}$ as
\[
-\ad_{\omega} \gamma = [\omega, \gamma] 
= \omega \gamma-\gamma\omega, \quad
\omega\,,\gamma \in \mathfrak{g},
\]
and define its dual $\operatorname{ad}^*:\mathfrak{g}\times\mathfrak{g}^*\to\mathfrak{g}^*$ as
\[
\langle \ad_\omega^*\mu, \gamma \rangle_{\mathfrak{g}} 
= 
\langle \mu, \ad_\omega \gamma \rangle_{\mathfrak{g}}
\,,\quad
\mu \in \mathfrak{g}.
\]
\end{definition}

\begin{lemma}[\bfi{Reduced equation for geodesic problem}] 
The Lie algebra variable $\xi$ for the geodesic matching problem stated in Definition \ref{matching problem} satisfies
\begin{equation}
\label{reduced}
\dd{}{t}\mathcal{A}\xi
 + \ad_{\xi}^*\mathcal{A}\xi = 0,
\end{equation}
weakly, in the sense of the pairing $\langle\,\cdot\,, \,\cdot\, \rangle_{\mathfrak{g}}:\mathfrak{g}\times\mathfrak{g}^*\to\mathbb{R}$. 
\end{lemma}
\begin{proof}
For any $\gamma\in\mathfrak{g}$, we have, upon substituting equation (\ref{xi}),
\begin{eqnarray*}
\iprod{\gamma}{\dd{}{t}(\mathcal{A}\xi)}{\mathfrak{g}}
& = & \dd{}{t}
\iprod{\gamma}{\mathbf{J}_G(p_q)}{\mathfrak{g}}
\\
& = & \dd{}{t}
\iprod{p_q}{\mathcal{X}^G_{\gamma}q}{T^*Q} \\
& = & \iprod{\dd{p_q}{t}}{\mathcal{X}^G_{\gamma}q}{T^*Q} 
+ \iprod{p_q}{T_{q}\left(\mathcal{X}^G_{\gamma}q\right)
\cdot\dd{q}{t}}{T^*Q} 
\\
& = & \iprod{-\left(T_{q}\left(\mathcal{X}^G_{\xi}q\right)\right)^*\cdot p_q}{\mathcal{X}_{\gamma}q}{T^*Q} 
+ \iprod{p_q}{T_{q}\left(\mathcal{X}^G_{\gamma}q\right)\cdot\mathcal{X}^G_{\xi}q
}{T^*Q} 
\\
& = & 
 \iprod{p_q}{-T_{q}\left(\mathcal{X}^G_{\xi}q\right)\cdot
\mathcal{X}^G_{\gamma}q
+T_{q}\left(\mathcal{X}^G_{\gamma}q\right)\cdot
\mathcal{X}^G_{\xi}q
}{T^*Q} \\
& = & \iprod{p_q}{\mathcal{X}_{[\gamma,\xi]}q}{T^*Q} \\
& = & \iprod{[\gamma,\xi]}{\mathbf{J}_H(p_q)}{\mathfrak{g}} \\
& = & \iprod{[\gamma,\xi]}{\mathcal{A}\xi}{\mathfrak{g}} \\
& = & -\iprod{\gamma}{\ad_{\xi}^*\mathcal{A}\xi}{\mathfrak{g}}.
\end{eqnarray*}
Consequently, we obtain the result stated, since $\gamma$ is an arbitrary element of
$\mathfrak{g}$.
\end{proof}

We will next define a modification of the problem stated in
Definition \ref{matching problem}, which has the advantage that the
endpoint conditions do not contain a free reparameterisation
variable. This reformulation is more amenable when solving the curve
matching problem numerically, for example. We shall proceed to show that 
solutions of the modified problem can be transformed into solutions of
the problem stated in Definition \ref{matching problem}.

\begin{definition}[\bfi{Reparameterised geodesic problem with
    symmetry}]
\label{reparameterised matching problem}
Let $Q$ be a manifold, let $G$ be a Lie group acting on $Q$ from the
left, and let $H$ be a (possibly different) Lie group acting on $Q$
from the right that commutes with the left action of $G$, with
corresponding Lie algebras $\mathfrak{g}$ and $\mathfrak{h}$, and
corresponding Lie algebra actions $\mathcal{X}^G$ and $\mathcal{X}^H$
respectively.  Furthermore, let
$\mathcal{A}:\mathfrak{g}\to\mathfrak{g}^*$ be a positive-definite
self-adjoint operator. We seek
\begin{itemize}
\item a one parameter family ${{q}}$ of points on $Q$ for $t\in[0,1]$,
\item a one parameter family ${{\xi}}$ of elements of
  $\mathfrak{g}$ for $t\in[0,1]$, and
\item $\nu\in \mathfrak{h}$, 
\end{itemize}
which minimise
\[
\int_0^1\frac{1}{2}\langle{{\xi}},\mathcal{A}{{\xi}}\rangle
_{\mathfrak{g}}\diff{t}
\,,
\]
where $\langle \cdot,\cdot\rangle_{\mathfrak{g}}$ is the usual
inner product on $\mathfrak{g}$ , subject to the constraints
\begin{eqnarray} \mbox{\rm [Reconstruction relation]} \quad
  \dd{}{t}{{q}}
  &=& \mathcal{X}^G_{{\xi}}{q} + \mathcal{X}^H_{\nu}{q}, 
\label{reparameterised dynamical constraint}  \\
  \mbox{\rm [Initial state (Template)]} \quad
  {{q}}_{t=0} & = & {q}^A, \label{reparameterised template constraint} \\
  \mbox{\rm [Final state (Target)]} \quad {{q}}|_{t=1} & = &
  {q}^B, \label{reparameterised target constraint}
\end{eqnarray}
and ${q}^A$, ${q}^B$ are chosen points on $Q$.
\end{definition}
\begin{remark}
  Note that in this modified definition, we do specify the final
  boundary condition for ${q}$ without allowing arbitrary symmetry
  transformations using $H$. However, we also introduce an additional
  variable $\nu$ which moves ${q}$ in the direction of symmetries
  generated by $\mathfrak{h}$.
\end{remark}
We shall derive the equations of motion associated with this modified
problem, and the associated conservation laws. These will lead us to
conclude that it possible to construct solutions of the problem in
Definition \ref{matching problem} out of solutions of the problem in
Definition \ref{reparameterised matching problem}, and the latter can be solved
as a shooting problem in which the boundary conditions are explicitly
specified, rather than as an algebraic condition. As before, we can
derive the equations of motion for $\bar{q}$, ${{\xi}}$ and the condition
for $\nu$ by enforcing the reconstruction relation
\eqref{reparameterised dynamical constraint} as a constraint using
Lagrange multipliers ${\bar{p}_q}\in T_{\bar{q}}^*Q$, leading to the
following variational principle.

\begin{definition}[\bfi{Variational principle for reparameterised
    geodesic problem with symmetry}]
\label{reparameterised variational principle matching problem}
We seek $({p},{q})\in T^*Q$ and ${\xi}\in\mathfrak{g}$ for
$t\in[0,1]$, and $\nu\in \mathfrak{h}$, which satisfy
\begin{eqnarray}
\delta S = \delta \int_0^1
\frac{1}{2}\left\langle {\xi}, \mathcal{A}{\xi}
\right\rangle_{\mathfrak{h}}
+ \left\langle {p}_q,\dd{}{t}{q} + \mathcal{X}^G_{{\xi}}{q} 
+ \mathcal{X}^H_{\nu}{q}\right\rangle
_{T^*Q}
\diff{t} = 0,
\end{eqnarray}
subject to 
\begin{eqnarray}
q|_{t=0} = q^A, \quad q|_{t=1} = q^B,
\end{eqnarray}
under variations of ${p}_q$, ${q}$, ${\xi}$ and $\nu$.
\end{definition}
Proceeding just as before, we can use variational calculus to obtain
the equations of motion, as described in the following lemma.
\begin{lemma}[\bfi{Equations of motion for reparameterised geodesic problem}]
\label{reparameterised matching problem eqns}
At the optimum, the following equations are satisfied in the sense of appropriate pairings:
\begin{eqnarray}
\label{reparameterised dot q}
\dd{}{t}{q} - \mathcal{X}^G_{{\xi}} {q} - \mathcal{X}^H_{\nu} {q}
& = & 0, \\
\label{reparameterised dot p}
\dd{}{t}{p}_q + T_{{q}}\left(\mathcal{X}_{{\xi}}^G {q} + \mathcal{X}_{\nu}^H{q}\right)^*
\cdot {p}_q
& = & 0, \\
\label{reparameterised xi}
\mathcal{A}{\xi} - \mathbf{J}_G({p}_q) & = & 0.
\end{eqnarray}
Furthermore,
\begin{equation}
\label{reparameterised nu}
\int_{t=0}^1
\mathbf{J}_H({p}_q)
\diff{t} = 0.
\end{equation}
\end{lemma}
\begin{proof}
\begin{eqnarray*}
\delta S 
%& = & 
%\int_0^1 \iprod{\delta {\xi}}{\mathcal{A}{\xi}}{\mathfrak{g}}
%+\iprod{\delta {p}_q}{\dd{{q}}{t}-\left(\mathcal{X}^G_{{\xi}}
%+ \mathcal{X}^H_{\nu}\right){q}}{T^*Q}
%+\iprod{{p}_q}{\delta\dd{{q}}{t}-\delta(\mathcal{X}^G_{{\xi}}\%overline{q})
%-\delta(\mathcal{X}^H_\nu {q})}{T^*Q}\diff{t} \\
& = & \int_0^1 \iprod{\delta {\xi}}{\mathcal{A}{\xi}}{\mathfrak{g}}
+\iprod{\delta {p}_q}{\dd{{q}}{t}-
\left(\mathcal{X}^G_{{\xi}}+\mathcal{X}^H_{\nu}\right){q}}{T^*Q}
+\iprod{{p}_q}{\delta\dd{{q}}{t}
-\left(\mathcal{X}^G_{\delta{\xi}}+\mathcal{X}^H_{\delta\nu}\right){q}
-T_{{q}}\left(\mathcal{X}^G_{{\xi}}{q} +
 \mathcal{X}^H_{\nu}{q}\right)\delta {q}}{T^*Q}
\diff{t} \\
& = & \int_0^1 \iprod{\delta {\xi}}{\mathcal{A}{\xi}
  - \mathbf{J}_G({p}_q)}
{\mathfrak{g}}
+\iprod{\delta {p}}{\dd{{q}}{t}-
\left(\mathcal{X}^G_{{\xi}}+\mathcal{X}^H_{\nu}\right){q}}{T^*Q}
+\iprod{-\dd{{p}}{t} -\left(T_{{q}}\left(\mathcal{X}^G_{{\xi}}{q} +
 \mathcal{X}^H_{\nu}{q}\right)\right)^*{p}}
{\delta {q}}{T^*Q}\diff{t} \\
& & \qquad + \iprod{\delta\nu}{\int_0^1 
\mathbf{J}_H({p}_{{q}})\diff{t}}{\mathfrak{h}}.
\end{eqnarray*}
Since $\delta {p}$, $\delta {q}$ and $\delta {\xi}$ are all arbitrary
we obtain equations (\ref{reparameterised dot q}-\ref{reparameterised
  nu}).
\end{proof}
Proceeding as before, we can transform (\ref{reparameterised nu})
into an initial condition by making use of the conservation of the
right-action momentum map, $\mathbf{J}_H$.
\begin{lemma}[\bfi{Noether's theorem for reparameterised geodesic problem}]
\label{reparameterised noether matching}
The system of equations (\ref{reparameterised dot
  p}-\ref{reparameterised xi}) has a conserved  momentum $\mathbf{J}_H$.
\end{lemma}
\begin{proof}
  The problem in Definition \ref{reparameterised variational principle matching
    problem} is invariant under transformations
\[
{q} \to {q}\alpha, \qquad \alpha\in H,
\]
which are generated by
$\gamma=\in\mathfrak{h}$. This
means that the variational principle in Definition
\ref{reparameterised variational
  principle matching problem} is invariant under application of the
\bfi{cotangent lift} (\emph{i.e.}, the dual of the inverse of its
infinitesimal transformation in $Q$) namely
\[
\delta {q} = \mathcal{X}^H_{\gamma}{q}, \quad \delta {p} =
-\left(T_{{q}}\left( \mathcal{X}^H_{\gamma}{q}\right)\right)^*{p},
\quad \delta\nu =0,
\]
where for convenience we have inserted the time dependence
\[
\gamma = 0 \quad\mbox{if}\quad t_0<t<1
.\]
 Substitution of this
infinitesimal transformation into the variational principle gives
\begin{eqnarray*}
0 = \delta S & = & 
\int_0^{t_0} 
\left[ \iprod{\delta {p},\dd{}{t}{q} - \mathcal{X}^G_{{\xi}}{q} - \mathcal{X}^H_{\nu}{q}}{T^*Q}
+ \iprod{{p}}{\delta \dd{}{t}{q} - 
\left(T_{{q}}\left(\mathcal{X}^G_{{\xi}}{q}\right) +
T_{{q}}\left(\mathcal{X}^G_{{\xi}}{q}\right)\right)
\cdot\delta {q}}{T^*Q}
\right]\diff{t}
\\
& = & \int_0^{t_0}
\left[
\iprod{\delta {p}}{\underbrace{\dd{}{t}{q} - \mathcal{X}^G_{{\xi}}{q} - \mathcal{X}^H_{\nu}{q}}_{\hbox{=0}}}{T^*Q}
-\iprod{\underbrace{\dd{}{t}{p} + 
\left(T_{{q}}\left(\mathcal{X}^G_{{\xi}}
+\mathcal{X}^H_{\nu}\right)
\right)^*{p}}_{\hbox{=0}}}{\delta {q}}{T^*Q}
\right]\diff{t} + \left[
\iprod{{p}}{\delta {q}}{T^*Q}
\right]_{t=0}^{t=t_0} \\
&=& \left[
\iprod{{p}}{\mathcal{X}^H_{\gamma}{q}}{T^*Q}
\right]_{t=0}^{t=t_0} \\
&=& \left[
\iprod{\gamma}{\mathbf{J}_H(p_q)}{\mathfrak{h}}\right]_{t=0}^{t=t_0}. \\
\end{eqnarray*}
Since this equation holds for any $0<t_0<1$, then the right-action
momentum map $\mathbf{J}_H$ is conserved.
\end{proof}

\begin{corollary}[\bfi{Vanishing momentum}]
\label{reparameterised vanishing momentum}
The conserved momentum satisfies $\mathbf{J}_H=0$ for
all times $t$.
\end{corollary}
\begin{proof}
  After noting that $\mathbf{J}_H$ is conserved, Lemma \ref{reparameterised matching problem eqns} gives
\[
0 = \int_0^1 \mathbf{J}_H dt = \mathbf{J}_H.
\]
\end{proof}
Next we show that $\xi$ obtained from Definition \ref{reparameterised
  variational principle matching problem} satisfies the same
Euler-Poincar\'e equation as $\xi$ obtained from Definition
\ref{variational principle matching problem}.
\begin{lemma}[\bfi{Reduced equation for geodesic problem}]
The Lie algebra variable ${\xi}$ obeys equation
\eqref{reduced} weakly, i.e., for an appropriate pairing.
\end{lemma}
\begin{proof}
For any $\gamma\in\mathfrak{g}$, we have, upon substituting equation (\ref{reparameterised xi}),
\begin{eqnarray*}
\iprod{\gamma}{\dd{}{t}(\mathcal{A}{\xi})}{\mathfrak{g}}
& = & \,\dd{}{t}\iprod{\gamma}{\mathbf{J}_G(p_q)}{\mathfrak{g}}
\\
& = & \dd{}{t}
\iprod{{p_q}}{\mathcal{X}^G_{\gamma}{q}}{T^*Q} \\
& = & \iprod{\dd{{p_q}}{t}}{\mathcal{X}^G_{\gamma}{q}}{T^*Q} 
+ \iprod{{p_q}}{T_{{q}}\left(\mathcal{X}^G_{\gamma}{q}\right)
\cdot\dd{{q}}{t}}{T^*Q} 
\\
& = & \iprod{-\left(T_{{q}}\left(\mathcal{X}^G_{{\xi}}{q}
+ \mathcal{X}^H_{\nu}{q}
\right)\right)^*\cdot {p}_q}{\mathcal{X}^G_{\gamma}{q}}{T^*Q} 
+ \iprod{{p_q}}{T_{{q}}\left(\mathcal{X}^G_{\gamma}{q}\right)
\cdot\left(\mathcal{X}^G_{{\xi}}{q}+\mathcal{X}^H_{\nu}{q}
\right)
}{T^*Q} 
\\
& = & 
\iprod{{p_q}}
{-T_{{q}}\left(\mathcal{X}^G_{{\xi}}{q}
+\mathcal{X}^H_\nu{q}\right)
\cdot
\mathcal{X}^G_{\gamma}{q}
+T_{{q}}\left(\mathcal{X}^G_{\gamma}{q}\right)\cdot
\left(\mathcal{X}^G_{{\xi}}{q}
+\mathcal{X}^H_\nu{q}\right)
}{T^*Q} \\
& = & \iprod{{p_q}}{\mathcal{X}^G_{[\gamma,{\xi}]}{q}}{T^*Q}
+ \underbrace{
\iprod{{p_q}}{T_{{q}}\left(\mathcal{X}^G_{\gamma}{q}\right)\cdot
\mathcal{X}^H_\nu{q}
- T_{{q}}\left(\mathcal{X}^H_{\nu}{q}\right)\cdot
\mathcal{X}^G_\gamma{q}
}{T^*Q}}_{\hbox{=0}} \\
& = & \iprod{[\gamma,{\xi}]}{\mathbf{J}_G(p_q)}{\mathfrak{h}} \\
& = & \iprod{[\gamma,{\xi}]}{\mathcal{A}{\xi}}{\mathfrak{h}} \\
& = & -\iprod{\gamma}{\ad_{{\xi}}^*\mathcal{A}{\xi}}{\mathfrak{h}},
\end{eqnarray*}
and we obtain the stated result, since $\gamma$ is an arbitrary element of
$\mathfrak{g}$. Here, the underbraced term vanishes since the left and right
group actions commute with each other.
\end{proof}
This means that we can show that the two problems produce equivalent
solutions provided that the initial conditions for $\xi$ are the same
in both cases. The following theorem establishes this result.
\begin{theorem}
\label{reparameterised theorem}
  Let $\overline{q}$, $\overline{p}_{\overline{q}}$, $\nu$, $\overline{\xi}$ be obtained from
  the solution of equations (\ref{reparameterised dot
    q}-\ref{reparameterised nu}), and define 
  \[
  \psi = \exp(\nu t) \in \mathcal{H}.
  \]
  Then the transformed variables constructed from
  \begin{equation}
    \label{cotangent lift}
  q = R_{\psi^{-1}}\overline{q}, \quad
  p_q = T^*_{q}\psi\left(\overline{p}_{\overline{q}}\right), \quad
  \xi = \overline{\xi}, \qquad t\in[0,1],
  \end{equation}
  (\emph{i.e.} the cotangent lift of $\psi$) satisfy equations
  (\ref{dot q}-\ref{xi}) together with the boundary conditions
  (\ref{template constraint},\ref{target constraint},\ref{end
    condition}), with $\eta=\psi^{-1}_{t=1}$.  Hence, $q$, $\eta$
  and $\xi$ form a (local) extremum for the problem in Definition
  \ref{variational principle matching problem}.
\end{theorem}
\begin{proof}
  First we take $\xi$ and $\overline{\xi}$ from equations \eqref{xi} and
  \eqref{reparameterised xi} respectively, and show that $\xi=\overline{\xi}$.
  Since the left and right actions commute, we have that
\begin{eqnarray*}
\iprod{\gamma}{\mathcal{A}\overline{\xi}}{\mathfrak{g}}
&=&\iprod{\overline{p}_{\overline{q}}}{\mathcal{X}^G_\gamma \overline{q}}{T^*M} \\
&=&\iprod{\overline{p}_{\overline{q}}}{\mathcal{X}^G_\gamma R_{\psi}q}{T^*M} \\
& = & \iprod{p_q}{T_{q}R_{\psi}\mathcal{X}^G_{\gamma}q}{T^*M} \\
& = & \iprod{T_{q}^*R_{\psi}(\overline{p}_{\overline{q}})}{\mathcal{X}^G_{\gamma}q}{T^*M} \\
& = & \iprod{p_q}{\mathcal{X}^G_{\gamma}q}{T^*M} \\
& = & \iprod{\gamma}{\mathcal{A}\xi}{\mathfrak{g}},
\end{eqnarray*}
and hence $\xi=\overline{\xi}$. Next we verify the equations for $q$ and $p_q$. 
Taking the time derivative of $\overline{q}$, we have
\begin{eqnarray*}
\dd{}{t}\overline{q} & = & \dd{}{t}\left(R_{\psi}q\right) \\
& = & T_{q}R_{\psi}\cdot\dd{}{t}q + 
\mathcal{X}^H_\nu\left(R_{\psi}q\right),
\end{eqnarray*}
and so
\begin{eqnarray*}
\dd{}{t}q &=& \left(T_{q}R_{\psi}\right)^{-1}
\left(\dd{}{t}\overline{q} - \mathcal{X}^H_{\nu}\overline{q}\right) \\
& = & \left(T_{q}R_{\psi}\right)^{-1}
\mathcal{X}^G_{\overline{\xi}}\overline{q} \\
& = & \mathcal{X}^G_{\overline{\xi}}q \\
& = & \mathcal{X}^G_{\xi}q,
\end{eqnarray*}
as required. To check the time evolution equation for $p_q$, we take
the inner product with an arbitrary tangent vector $v\in T_{q}Q$, to find:
\begin{eqnarray*}
\dd{}{t}\iprod{p_q}{v}{T^*Q} & = & 
\dd{}{t}\iprod{T^*_{q}R_{\psi}(\overline{p}_{\overline{q}})}{v}{T^*Q} \\
& = & \dd{}{t}\iprod{\overline{p}_{\overline{q}}}{\left(T_{q}R_{\psi}\right)\cdot v}{T^*Q} \\
& = & \iprod{\dd{}{t}\overline{p}_{\overline{q}}}{\left(T_{q}R_{\psi}\right)\cdot v}{T^*Q} 
+ \iprod{\overline{p}_{\overline{q}}}{\dd{}{t}\left(T_{q}R_{\psi}\right)\cdot v}{T^*Q} \\
& = & -\iprod{T^*_{\overline{q}}\left(\mathcal{X}^G_{\overline{\xi}}\overline{q}
+ \mathcal{X}^H_{\nu}\overline{q}\right)\cdot\overline{p}_{\overline{q}}}{\left(T_{q}
R_{\psi}\right)\cdot v}{T^*Q} \\
& & \qquad +\iprod{\overline{p}_{\overline{q}}}{T_{\overline{q}}\left(\mathcal{X}^H_{\nu}\overline{q}\right)\cdot(T_{q}R_{\psi})\cdot v}{T^*Q} \\
& = & -\iprod{\overline{p}_{\overline{q}}}{T_{\overline{q}}\left(\mathcal{X}^G_{\overline{\xi}}\overline{q}\right)\cdot\left(T_{q}R_{\psi}\right)\cdot v}{T^*Q} \\
& = & -\iprod{T^*_{q}R_{\psi}(\overline{p}_{\overline{q}})}{T_{q}\left(\mathcal{X}^G_{\overline{\xi}}q\right)\cdot v}{T^*Q} \\
& = & -\iprod{p_q}{T_{q}\left(\mathcal{X}^G_{\overline{\xi}}q\right)\cdot v}{T^*Q} \\
& = & -\iprod{T^*\left(\mathcal{X}^G_{\overline{\xi}}q\right)\cdot p_q}{v}{T^*Q},
\end{eqnarray*}
as required. 

It remains to check the boundary conditions. Trivially, $q|_{t=0} =
q^A$, $q|_{t=1} = \overline{q}|_{t=1}\psi|_{t=1}^{-1} = q^B\eta$, as
required.  Finally, we need to check the end condition \eqref{end
  condition}. From Corollary \ref{reparameterised vanishing momentum},
we have
\[
\mathbf{J}_H(p_q)|_{t=0} =
\mathbf{J}_H(\overline{p}_{\overline{q}})|_{t=0} = 0,
\]
and Lemma \ref{noether matching} implies that $\mathbf{J}_H(p_q)|_{t=1}=0$.  Hence the boundary conditions are satisfied.
\end{proof}

\section{Examples}
\label{examples}
In this section we describe examples of the reparameterised geodesic
problem with symmetry, and discuss its applications to the
characterisation of the shape of curves and surfaces. 

\subsection{Example: SO(3)}
We illustrate our results with the case of the action of $SO(3)$ on
itself which gives rise to the equations of a rotating rigid body. We
consider the problem in which the end point boundary condition is only
determined up to a rotation of the rigid body about its $z$-axis. Of
course, this problem can also be solved by picking reduced
coordinates, but we use it as here as a simple example.
\begin{definition}[\bfi{Optimal control of a symmetric rigid body}]
 \label{egg rotating}
 Let $Q(t)$ be a one-parameter family of matrices in $SO(3)$.  Let
 $\omega(t)$ be a one-parameter family of matrices in
 $\mathfrak{so}(3)$. Let $R_\theta$ be a rotation in the $z$-axis by
 an angle $\theta$. We seek $Q(t)$, $\omega(t)$, and $R_\theta$ which
 satisfy
 \[
 \min_{\omega,\theta}\int_0^1
 \frac{1}{2}\Big\langle
 \omega,I\omega\Big\rangle
 _{\mathfrak{so}(3)\times\mathfrak{so}(3)^*}
 \diff{t}
 \]
 subject to the constraints
\begin{eqnarray} \mbox{\rm [Reconstruction relation]} \quad 
 \dot{Q}(t)
   &=& \omega(t)Q(t), \label{dynamical constraint egg}  \\
 \mbox{\rm [Initial (Template)]} \quad 
 Q(0) & = & Q^0, \label{template constraint egg} \\
 \mbox{\rm [Final (Target)]} \quad 
 Q(1) & = & Q^1R_\theta, \label{target constraint egg}
 \end{eqnarray}
 where $I$ is a chosen symmetric matrix. The dynamical constraint
 (\ref{dynamical constraint egg}) allows the reconstruction of the
 curve $Q(t)\in SO(3)$ on the Lie group from the optimal
 right-invariant (spatial) angular frequency
 \[
 \omega(t)=\dot{Q}Q^{-1}(t)
 \,,
 \]
 in the Lie algebra $\mathfrak{so}(3)$. The other constraints specify
 the starting and ending points of the curve $Q(t)\in SO(3)$.
\end{definition}
This problem is an example of the optimal control problem in
Definition \ref{matching problem}, with the manifold $Q$ being
$SO(3)$, the group $G$ being $SO(3)$ acting from the left, and the
group $H$ being $SO(2)$ acting from the right.  We identify
\[
q\equiv Q, \quad \omega\equiv \xi, \quad R_\theta \equiv \eta, \quad 
\mathcal{A}\equiv I
\mbox{ and } p_q\equiv P,
\]
where $P\in T^*SO(3)$ is the conjugate momentum to $Q$. Application of
Lemma \ref{matching problem eqns} gives the following dynamical
equations:
 \begin{eqnarray}
 \label{dot Q}
 \dot{Q} - \omega Q & = & 0, \\
 \label{dot P}
 \dot{P} + \omega^T P & = & 0, \\
 \label{omega}
 I\omega + \frac12(PQ^T-QP^T) & = & 0,
% P_{ij}w_{ik}Q_{kj} = w_{ij}P_{ik}Q_{jk}
 \end{eqnarray}
which are the equations for a rotating rigid body.
The last line follows from the definition of the left-action momentum map for $SO(3)$, namely
\[
\iprod{\mathbf{J}^L(P)}{\delta \omega}{\mathfrak{so}(3)}
=
\iprod{P }{-{\delta \omega} Q}{T^*SO(3)}
=
-{\rm tr}\Big(P^T{\delta \omega} Q\Big)
=
-{\rm tr}\Big((PQ^T)^T{\delta \omega} \Big)
=
\iprod{PQ^T }{-{\delta \omega}}{\mathfrak{so}(3)}
\,,
\]
where $\delta \omega$ is an arbitrary antisymmetric $3\times3$ matrix.
(Hence, the antisymmetric combination in equation (\ref{omega}).)
The end point condition (which comes from minimising over $\theta$) becomes
\begin{equation}
\label{SO(3) end condition}
\iprod{P}{Qw}{T^*SO(3)}|_{t=1}=0,
\end{equation}
where
\begin{equation}
\label{w definition}
w = \begin{pmatrix}
0 & 1 & 0 \\
-1 & 0 & 0 \\
0 & 0 & 0 \\
\end{pmatrix},
\end{equation}
% w_{ij} = E_{ijk}n_k, n = (0,0,1)^T
and Lemma \ref{vanishing momentum} implies that the quantity 
$\iprod{P}{Qw}{T^*SO(3)}$ (which is the $z$-component of the angular momentum) is zero for all times $t$. 
Lemma \ref{reduced}
states that $\omega$ satisfies the Euler equations for a rigid body:
\[
\dd{}{t}(I\omega) + \ad_\omega^*(I\omega) = 0,
\]
where we define
\[
\ad_{\omega} \gamma = [\omega, \gamma] 
= \omega \gamma-\gamma\omega, \quad
\omega\,,\gamma \in \mathfrak{so}(3),
\]
and
\[
\langle \ad_\omega^*\mu, \gamma \rangle = \langle \mu, \ad_\omega \gamma
\rangle
\,,\quad
\mu \in \mathfrak{so}(3)^*.
\]

For this problem, obtaining a solution is simple, since one can define
coordinates on $T(SO(3))$, and remove the coordinates associated with
the $R_\theta$ direction and the corresponding vanishing conserved
momentum, and solve a two-part boundary problem for the remaining
coordinates. However, we wish to develop a methodology for numerical
discretisations of infinite-dimensional problems where it is less
clear how to do this. Hence, we define the following formulation which
makes use of a time-varying ``relabelling'' transformation in the
$R_\theta$ direction. Theorem \ref{reparameterised theorem} states
that to obtain solutions to equations (\ref{dot Q}-\ref{omega}), we
can solve the following modified problem:
\begin{definition}[\bfi{Reparameterised optimal control of a symmetric
    rigid body}]
 \label{reparameterised egg rotating}
   Let $\overline{Q}(t)$ be a one-parameter family of matrices in $SO(3)$.
   Let $\overline{\omega}(t)$ be a one-parameter family of matrices in
   $\mathfrak{so}(3)$. Let $\nu$ be the generator of a rotation
   about the $z$-axis, which may be written in the form
   \[
   \nu = \theta w,
   \]
   where $w$ is defined in equation \eqref{w definition}, and
   $\theta\in\mathbb{R}$.

   We seek $\overline{Q}(t)$, $\overline{\omega}(t)$, and $\nu$ which satisfy
   \[
   \min_{\overline{\omega},\nu}\int_0^1
   \frac{1}{2}\Big\langle
   \overline{\omega},I\overline{\omega}\Big\rangle
   _{\mathfrak{so}(3)\times\mathfrak{so}(3)^*}
   \diff{t}
 \]
 subject to the constraints
\begin{eqnarray} \mbox{\rm [Reconstruction relation]} \quad 
 \dot{\overline{Q}}(t)
   &=& \overline{\omega}(t)\overline{Q}(t) + \overline{Q}(t)\nu, \label{reparam dynamical constraint egg}  \\
 \mbox{\rm [Initial (Template)]} \quad 
 \overline{Q}(0) & = & Q^0, \label{reparam template constraint egg} \\
 \mbox{\rm [Final (Target)]} \quad 
 \overline{Q}(1) & = & Q^1, \label{reparam target constraint egg}
 \end{eqnarray}
 where $I$ is a chosen symmetric matrix.
\end{definition}
Lemma \ref{reparameterised matching problem eqns} states that the solution
to this problem satisfies the following equations:
\begin{eqnarray}
\label{dot Qoverline}
\dot{\overline{Q}} - \overline{\omega} \overline{Q} + \overline{Q}\nu & = & 0, \\
\label{dot Poverline}
\dot{\overline{P}} + \overline{\omega}^T \overline{P} - \overline{P}\nu^T & = & 0, \\
\label{omega relabelled}
I\overline{\omega} + \frac12(\overline{P}\,\overline{Q}^T-\overline{Q}\,\overline{P}^T) & = & 0,
\end{eqnarray}
with end-point condition
\begin{equation}
\iprod{\overline{P}}{\overline{Q}w}{T^*SO(3)}|_{t=1}=0.
\end{equation}
Lemma \ref{reparameterised vanishing momentum} states that 
\begin{equation}
\label{SO(3) vanishing momentum}
\iprod{\overline{P}}{\overline{Q}w}{T^*SO(3)}=0,
\end{equation}
for all $t$.

Hence, to obtain a solution to equations (\ref{dot Q}-\ref{SO(3) end
  condition}), we solve the two-point boundary value problem given by
equations (\ref{dot Qoverline}-\ref{dot Poverline}) with boundary
conditions (\ref{reparam template constraint egg}-\ref{reparam target
  constraint egg}). This can be formulated as a shooting problem, in
which we seek $\nu$ (or, equivalently, $\theta$) and initial
conditions for $P$ satisfying equation \eqref{SO(3)
  vanishing momentum}, such that the end point boundary
condition \eqref{reparam target constraint egg} is satisfied. We then
construct the reparameterisation matrix $R(t)$ from
\[
R(t) = \exp(\nu t),
\]
and use equation \eqref{cotangent lift} to reconstruct the solution in the form:
\[
Q(t) = \overline{Q}(t)R(t)
\quad \hbox{and}\quad
P(t) = \overline{Q}(t)R^T(t), 
\quad \hbox{implying}\quad
\omega(t) = \overline{\omega}(t)
\,,
\]
since, e.g., $R^T=R^{-1}$ implies $\overline{P}\,\overline{Q}^T=PQ^T$.
Then substituting these relations into equations (\ref{dot Q}-\ref{omega}) and equation \eqref{SO(3) end condition} recovers equations (\ref{dot Qoverline}-\ref{omega relabelled}) and equation (\ref{SO(3) vanishing momentum}). 

\subsection{Example: SE(3)}
We next describe the example of the action of $SE(3)$ on itself from
the left, with $SO(2)$ acting from the right. This example could describe a
docking problem of a spacecraft onto a space station. The spacecraft
can apply torque to rotate about a central point, or can apply thrust
to move itself in the direction in which it is pointing, and we wish
to dock the spacecraft using minimal energy. In the language of image
registration, this is known as rigid registration. We consider the
problem in which the end point boundary condition is only determined
up to a rotation of the rigid body about its $z$-axis. In the
spacecraft analogy, this corresponds to a docking procedure which does
not require the spacecraft to have any particular orientation about
the $z$-axis when docking. As in the previous example, this problem
can also be solved by picking reduced coordinates, but it serves as a
prototype for infinite dimensional problems.

Following the notation of \cite{Ho2009}, we represent an element $q$ of
$SE(3)$ as a $4\times4$ matrix:
\[
q \equiv \begin{pmatrix}
Q & \MM{r} \\
0 & 1 \\
\end{pmatrix},
\]
where $Q$ is an orthogonal matrix, and $\MM{r}\in \mathbb{R}^3$.
We represent an element $p_q$ of $T^*_qSE(3)$ as a $4\times 4$ matrix:
\[
p_q \equiv \begin{pmatrix}
P & \MM{p} \\
0 & 0 \\
\end{pmatrix},
\]
where $P\in T^*SO(3)$ and $\MM{p}\in\mathbb{R}^3$.  Finally, we
represent an element of the corresponding Lie algebra
$\mathfrak{se}(3)$ as another $4\times4$ matrix:
\[
\xi \equiv \begin{pmatrix}
\omega & \MM{v} \\
0 & 0 \\
\end{pmatrix},
\]
where $\omega$ is an antisymmetric matrix, and $\MM{v}\in
\mathbb{R}^3$. The reconstruction relation is then given by
\[
\dot{q} = \mathcal{X}^G_{\xi}q = \xi q = 
\begin{pmatrix}
\omega Q & \omega \MM{r} + \MM{v} \\
0 & 0 \\
\end{pmatrix}.
\]
We write the energy cost for the system as
\[
S = \int_{t=0}^1{E}\diff{t} = 
\frac12\int_{t=0}^1\iprod{\xi}{A\xi}{\mathfrak{se}(3)}\diff{t},
\]
where $A$ is a symmetric positive definite $4\times 4$ matrix given by
\[
A = 
\begin{pmatrix}
I & \mathbf{b} \\
\mathbf{b}^T & c \\
\end{pmatrix}.
\]
The
starting condition is specified as
\[
q|_{t=0} = 
\begin{pmatrix}
Q^A & \MM{r}^A \\
0 & 1 \\
\end{pmatrix},
\]
which specifies the starting orientation and position, and the end
condition is specified as
\[
q|_{t=1} = 
\begin{pmatrix}
Q^BR_{\theta} & \MM{r}^B \\
0 & 1 \\
\end{pmatrix},
\]
where $R_{\theta}$ is a rotation through any angle $\theta$ about the
$z$-axis. If we solve the problem in Definition \ref{matching
  problem}, then Lemma \ref{matching problem eqns} gives the 
dynamical equations
\begin{eqnarray}
\dot{Q} & = & \omega Q, \label{dot Q SE(3)} \\
\dot{\MM{r}} & = & \omega\MM{r} + \MM{v}, \\
\dot{P} & = & -\omega^T P, \\
\dot{\MM{p}} & = & \omega^T\MM{p}, \\
\frac{\delta E}{\delta \omega}
+ \frac12(PQ^T - QP^T)& = & 0, \\
\frac{\delta E}{\delta \mathbf{v}}
+ \MM{p} & = & 0, \label{v SE(3)}
\end{eqnarray}
where ${\delta E}/{\delta \omega}=
I{\omega} + \mathbf{b}\mathbf{v}^T -  \mathbf{v}\mathbf{b}^T$, 
${\delta E}/{\delta \mathbf{v}}=c\mathbf{v}+2\omega\mathbf{b}$.
The corresponding end point condition is
\begin{equation}\label{steve}
\iprod{P}{Qw}{T^*SO(3)}|_{t=1}=0,
\end{equation}
where $w$ is defined in equation \eqref{w definition}, as for the
$SO(3)$ case. From Lemma \ref{noether matching}, this quantity 
vanishes for all $t$. Theorem \ref{reparameterised theorem} then states
that a solution to these equations can be obtained by solving
the following reparameterised equations:
\begin{eqnarray}\label{colin}
\dot{\overline{Q}} & = & \overline{\omega} \overline{Q} + \alpha \overline{Q}w, \\
\dot{\overline{\MM{r}}} & = & \overline{\omega}\,\overline{\MM{r}} + \overline{\MM{v}}, \\
\dot{\overline{P}} & = & -\overline{\omega}^T \overline{P} - \alpha \overline{Q}w^T, \\
\dot{\overline{\MM{p}}} & = & \overline{\omega}^T\overline{\MM{p}}, \\
\dede{E}{\overline{\omega}} 
+ \frac12(\overline{P}\,\overline{Q}^T - \overline{Q}\,\overline{P}^T) & = & 0, \\
\dede{E}{\overline{\MM{v}}} + \overline{\MM{p}} & = & 0,
\label{rebecca}
\end{eqnarray}
with $\alpha\in \mathbb{R}$, $\overline{Q}|_{t=0}=Q^A$,
$\overline{\MM{r}}|_{t=0}=\MM{r}^A$, $\overline{Q}|_{t=1}=Q^B$,
$\overline{\MM{r}}|_{t=1} =\MM{r}^B$, and endpoint condition
\begin{equation}\label{allan}
\iprod{P}{Qw}{T^*SO(3)}|_{t=0}=0
\,.
\end{equation}
 This gives a two-point boundary
value problem with a constraint on the initial conditions plus an
extra parameter, which can be solved as a shooting problem by finding
$\alpha$ and $\overline{P}$ (subject to the constraint) such that $\overline{Q}$
and $\overline{\MM{r}}$ reach their target values $Q^B$ and $\MM{r}^B$.
A solution to the problem in Definition \ref{matching problem} can
then be reconstructed by defining $R(t)=\exp(\alpha w t)$, and
using the following formulae:
\[
Q(t) = \overline{Q}(t)R(t), \quad P(t) = \overline{P}(t)R^T(t), 
\quad \MM{r}(t)=\overline{\MM{r}}(t), \quad \MM{p}(t) = \overline{\MM{p}}(t), 
\quad \omega(t)=\overline{\omega}(t), \quad \MM{v}(t) = \overline{\MM{v}}(t). 
\]
Substituting these relations into equations (\ref{dot Q SE(3)}-\ref{v SE(3)}) and equation \eqref{steve} recovers equations (\ref{colin}-\ref{rebecca}) and equation (\ref{allan}). 

\subsection{Curve matching}\label{sec: curve matching}

In this section we return to the problem described in Definition
\ref{curve matching problem}, and discuss a number of practical issues
which are addressed by the   formulation discussed in
this paper. The aim of solving the problem is to find a
characterisation of the simple closed curve $\Gamma^B$ in terms of the
reference simple closed curve $\Gamma^A$, together with a scalar
periodic function $p(s)$ which specifies the initial conditions for
the normal component of the conjugate momentum $\MM{p}(s;t)$. In this
case, $Q$ is the space $\Emb(S^1,\mathbb{R}^2)$ of functions
$\MM{q}:S^1\to\mathbb{R}^2$, $G$ is the group $\Diff(\mathbb{R}^2)$ of
diffeomorphisms of $\mathbb{R}^2$ which acts on $Q$ from the left
\[
\Phi_L(g,\MM{q})(s) = g(\MM{q}(s)), \quad \forall s \in S^1,
\]
and $H$ is the group $\Diff(S^1)$ of diffeomorphisms of $S^1$, which
acts on $Q$ from the right
\[
\Phi_R(g,\MM{q})(s) = \MM{q}(g(s)), \quad \forall s \in S^1.
\]
The left and right actions can be expressed succinctly as 
\[
GQ = \Emb(S^1,  G\cdot \mathbb{R}^2)
\,,\quad
HQ=\Emb(H\cdot S^1,  \mathbb{R}^2).
\]
It is clear from this that the actions of $G$ and $H$ commute with
each other.

Lemma \ref{matching problem eqns} gives the dynamical equations
\begin{eqnarray*}
\pp{}{t}\MM{q}(s;t) & = & \MM{u}(\MM{q}(s;t),t), \\
\pp{}{t}\MM{p}(s;t) & = & -(\nabla\MM{u}(\MM{q}(s;t),t))^T\cdot\MM{p}(s;t)\,, 
\end{eqnarray*}
where the velocity $\MM{u}$ is defined weakly from the equation
\[
\iprod{\MM{w}}{A\MM{u}(\cdot,t)}{\mathfrak{X}} = \int_{S^1}\MM{p}(s;t)
\cdot\MM{w}(\MM{q}(s;t),t)\diff{s},
\]
where $\iprod{\cdot}{\cdot}{\mathfrak{X}}$ is the inner product
on the vector fields $\mathfrak{X}(\Omega)$
associated with the norm $\|\cdot\|_{\mathfrak{X}}$, for any test
function $\MM{w}\in \mathfrak{X}^*$. This equation has the weak
solution
\begin{equation}
m(x,t) = \sum  p(t,s)\delta(x-q(t,s))\,,
\label{sing-soln}
\end{equation}
which is the singular-solution momentum map, $\mathbf{J}_{\rm Sing}$
discussed in \cite{HoMa2004}.  The end condition is
\begin{equation}
\label{p end condition}
\MM{p}(s;1)\cdot\pp{}{s}\MM{q}(s;1)=0, \quad \forall s\in S^1,
\end{equation}
and Lemma \ref{noether matching} states that this conserved momentum
vanishes for all $t$. This corresponds to $\MM{p}(s;t)$ being normal
to the curve parameterised by $\MM{q}(s;t)$. Hence, to find geodesics
between $\Gamma^A$ and $\Gamma^B$, we solve a shooting problem and
seek initial conditions $\MM{p}(s;0)$ with vanishing tangential
component, which fix $\MM{q}(s;1)=\MM{q}^B(\eta(s))$ for some
$\eta\in\Diff(S^1)$. Having solved this problem, one can describe
$\Gamma^B$ entirely in terms of $p(s)=\MM{p}(s;0)\cdot\MM{n}(s)$ where
$\MM{n}$ is the normal to $\Gamma^A$. The solution to the problem
also provides the distance between the two curves.

There are a number of difficulties with solving such a shooting
problem numerically. The parameterisation of the curve must
necessarily be discretised numerically, typically by a list of points,
as in \cite{McMa2006,Co2008}, which can be obtained from a
piecewise-constant representation of $\MM{q}$ as a function of $s$
\cite{Vi2009}, or as piecewise linear geometric currents
\cite{VaGl2005}. Having taken the discretisation, the
reparameterisation symmetry is broken (although a remnant of it is
left behind, as described in \cite{Co2008}) which means that it is
difficult to obtain a discrete form of the end condition $\MM{q}(s;1)=
\MM{q}^B(\eta(s))$. As described in the Introduction, this problem has
been approached by proposing various functionals which are minimised
when $\MM{q}(s;1)$ and $\MM{q}(s)$ overlap the most. However, these
functionals produce quite a complicated landscape with local minima,
and the case of studying large deformations we have found that they
can result in odd artifacts in the shooting process. Also, the changes
in these functionals with respect to measurement error are quite
technical to quantify which makes statistical inference more
complicated.\\[1mm]

\begin{figure}
\begin{center}
\includegraphics*[width=7cm]{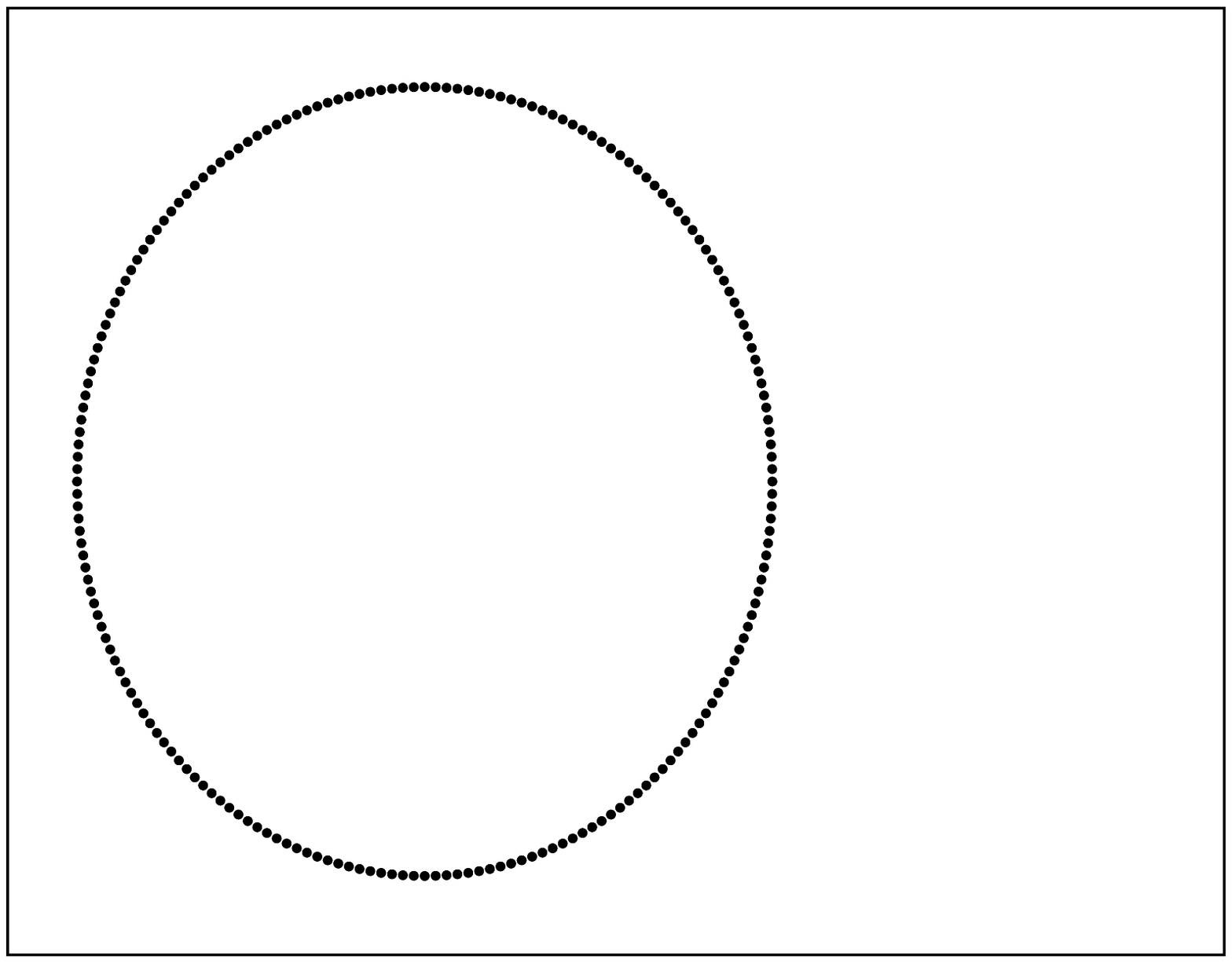}
\includegraphics*[width=7cm]{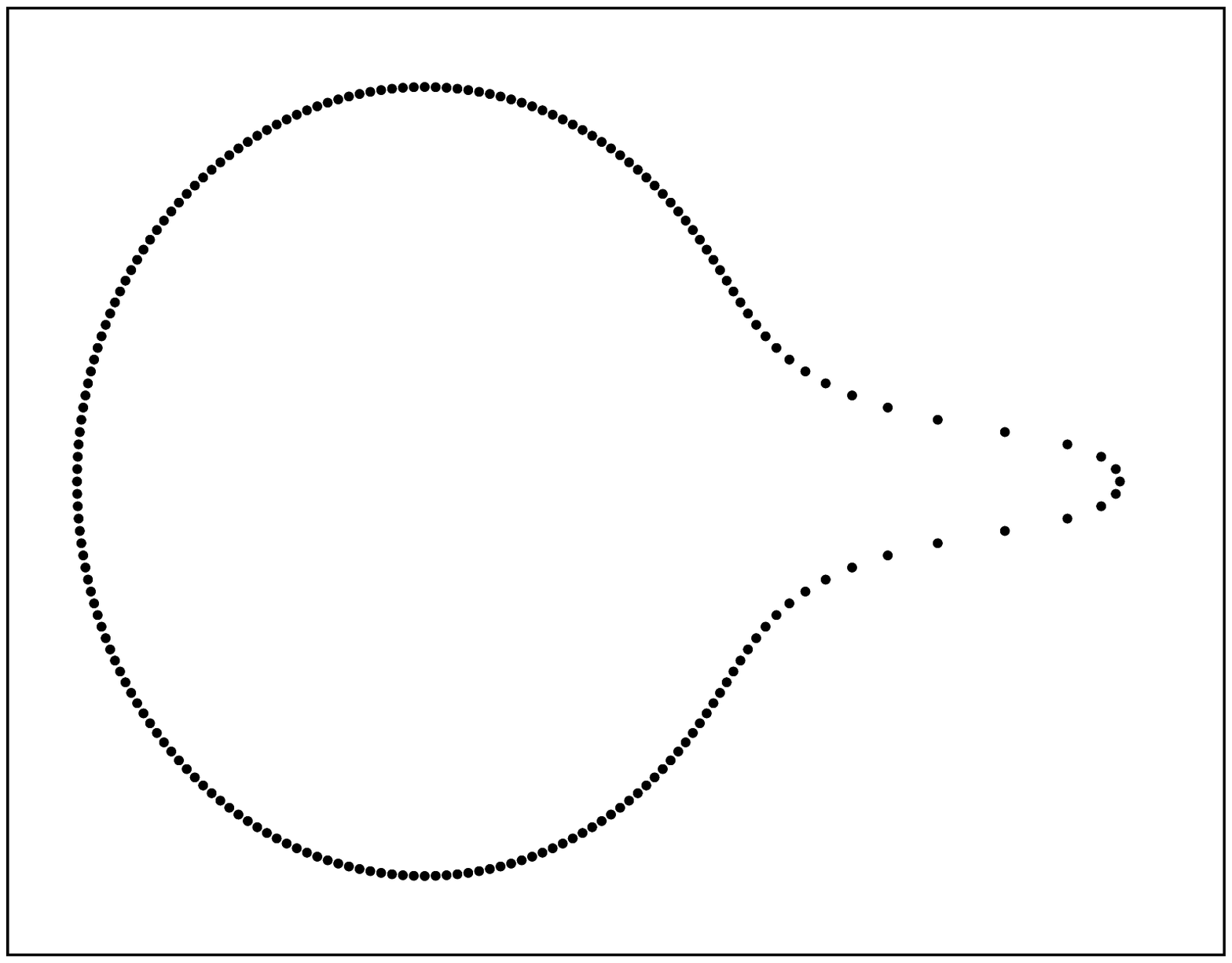}
\end{center}
\caption{\label{bulge} Figure illustrating the way in which
  deformation takes place in the parameterisation-independent geodesic
  equations for embedded curves. The initial curve is shown on the
  left, and the final curve is shown on the right.  Since the momentum
  is constrained to be normal to the curve, and since the velocity is
  obtained by applying a smoothing kernel to the momentum, the change
  in the shape emerges as local stretching of the curve, and the
  discrete points defining the shape become separated.}
\end{figure}

Another difficulty is that of adaptivity. As illustrated in Figure
\ref{bulge}, constraining $\MM{p}$ to be normal to the curve means
that any local large deformations give rise to large amounts of
stretching which then results in loss of accuracy in the approximation of
the functional used to enforce the end condition for $\MM{q}$. One
possible way to avoid this is to adaptively refine the grid point
density in the initial curve during the shooting process. \medskip

These difficulties are all removed if, instead, one solves the
following problem:
\begin{definition}[\bfi{Reparameterised curve matching problem}]
  \label{reparameterised curve matching problem}
  Let $\overline{\MM{q}}(s;t)$ be a one-parameter family of parameterised curves
  in $\mathbb{R}^2$, with $s\in [0,1]$ being the curve parameter and
  $t \in [0,1]$ being the parameter for the family. Let
  $\overline{\MM{u}}(\MM{x};t)$ be a one-parameter family of vector fields on
  $\mathbb{R}^2$. Let $\nu$ be a vector field on $S^1$. We seek
  $\overline{\MM{q}}$, $\overline{\MM{u}}$ and $\nu$ which satisfy
\[
\min_{\overline{\MM{u}},\nu}\int_0^1\frac{1}{2}\|\overline{\MM{u}}\|^2_{V}\diff{t}
\]
subject to the constraints
\begin{eqnarray} \mbox{\rm [Reconstruction relation]} \quad 
\pp{}{t}\overline{\MM{q}}(s;t)
  &=& \overline{\MM{u}}(\overline{\MM{q}}(s;t),t) + \nu(s)\pp{}{s}\MM{q}(s;t), \\
\mbox{\rm Initial state (Template)]} \quad 
\overline{\MM{q}}(s;0) & = & \overline{\MM{q}}^A(s), \\
\mbox{\rm [Final state (Target)]} \quad 
\overline{\MM{q}}(s;1) & = & \overline{\MM{q}}^B(s), \\
\end{eqnarray}
where $\|\cdot\|_{V}$ is the chosen norm which defines the space of
vector fields $V$.
\end{definition}
Lemmas \ref{reparameterised matching problem eqns} and
\ref{reparameterised vanishing momentum} state that extrema of this 
problem can be obtained by solving the shooting problem
\begin{eqnarray*}
  \pp{}{t}\overline{\MM{q}}(s;t) & = & \overline{\MM{u}}(\overline{\MM{q}}(s;t),t)
  + \nu(s)\pp{}{s}\overline{\MM{q}}(s;t), \\
  \pp{}{t}\overline{\MM{p}}(s;t) & = & -\overline{\MM{p}}(s;t)\cdot\nabla\overline{\MM{u}}(\overline{\MM{q}}(s;t),t) - \pp{}{s}\left(\nu(s)\overline{\MM{p}}(s;t)\right), \\
  \iprod{\MM{w}}{\overline{\MM{u}}(\cdot,t)}{V} &=& \int_{S^1}\overline{\MM{p}}(s;t)
  \cdot\MM{w}(\overline{\MM{q}}(s;t),t)\diff{s}, \qquad \forall \MM{w}\in V^*,
\end{eqnarray*}
with boundary conditions 
\[
\overline{\MM{q}}(s;0)=\MM{q}^A(s), \quad
\overline{\MM{p}}(s;0)\cdot\pp{}{s}\MM{q}^A = 0,
\quad \overline{\MM{q}}(s;1) = \MM{q}^B(s).
\]
The aim is to find $\nu(s)$ and normal components of $\MM{p}(s;0)$
such that these boundary conditions are satisfied. Note that in this
modified problem, the boundary condition for $\MM{q}(s;1)$ is
specified pointwise (\emph{i.e.} there is no reparameterisation
variable $\eta$ in the boundary condition). This means that the error
can be described directly in terms of 
\[
\|\MM{q}(\cdot;1)-\MM{q}^B\|^2
\]
for some chosen norm, which can be discussed in terms of measurement
errors directly. 

Theorem \ref{reparameterised theorem} then states that a solution to
the problem described in Definition \ref{curve matching problem} can
be reconstructed \emph{via}
\[
\MM{q}(s;t) = \overline{\MM{q}}(\eta(s;t),t), \quad
\MM{p}(s;t) = \overline{\MM{p}}(\eta(s;t),t)\pp{}{s}\eta(s;t)
\]
where
\[
\eta(s;t) = \exp(\nu(s)t).   
\]
This transformation produces a equivalent shooting problem in which the end condition for $\mathbf{\overline{q}}$ is now fixed.

\section{Summary and outlook}\label{summary}

In this paper, we studied an optimal control problem on a Lie group in
which the end boundary condition is specified only up to a symmetry.
We showed how this problem can be transformed into a modified problem
in which the end boundary condition is fixed, but an extra parameter
is introduced in the reconstruction relation, and proved that the two
problems are equivalent. This approach is motivated by the problem of
computing a geodesic on the diffeomorphism group in the plane which
takes one curve to another. The transformation gives rise to a system
of equations for a parameterised curve in which the end boundary
condition for each value of the parameter is fixed.  This means that
when a discrete approximation of the curve is used to solve this
problem numerically, the end boundary condition can again be fixed
exactly. In particular, when solving a shooting problem, this means
that the error between the computed curve and the target curve can be
computed simply by measuring the Euclidean distance between points for
each value. This method extends straightforwardly to the problem of
finding geodesics in the three-dimensional diffeomorphism group which
take one parameterised surface to another, with the end boundary
condition specified only up to reparameterisations of the target
surface. This problem has many applications in, for example,
biomedicine, since it allows topologically equivalent surfaces to be
characterised by a momentum field distributed on the template
surface. Such momentum fields exist in a linear space and so can be
manipulated using linear techniques
and still a topologically equivalent surface will always be obtained.

%Darryl:\\[2mm]
%(D1) framework left right dual / complementary dynamics: \\Both Left and Right %actions: left for motion and right for reparameterising.\\
%(D2) More broadly, this approach provides a general philosophy for image analys%is, e.g., medical images\\
%(D3) Potentially allows split between segmentation and registration\\

In the standard approach to planar image registration, the problem of
registering a specified closed curve (called the template) at time
$t=0$ to another (the target) at time $t=1$ amounts to deforming the
space $\mathbb{R}^2$ in which the template curve is embedded until it
matches the \emph{fixed} target image to within a certain
tolerance. Here, we considered a manifold $Q$ comprising the space of
closed curves $S^1$ embedded into the plane $\mathbb{R}^2$, written as
$Q=\Emb(S^1, \mathbb{R}^2)$. There are two Lie group actions available
for manipulating the closed curves in this description. The action
$G\times \mathbb{R}^2\to \mathbb{R}^2$ of the Lie group
$G=\Diff(\mathbb{R}^2)$ by composition from the \emph{left} deforms
the range space $\mathbb{R}^2$, and thereby drags along a curve
embedded in it as $GQ = \Emb(S^1, G\cdot \mathbb{R}^2)$. This left
action produced the singular-solution momentum map, $\mathbf{J}_{\rm
  Sing}$ in equation (\ref{sing-soln}), which introduces the
parameterisation of the closed curve by its position and momentum
supported on a delta function defining the curve in $\mathbb{R}^2$. Under this left
action of $G$, the curve preserves the initial parameterisation of its
domain space in $S^1$, although the current positions of the $S^1$
labels in the plane $\mathbb{R}^2$ will change as the range space is
transformed. Alternatively, the action $H\times S^1\to S^1$ of the Lie
group $H=\Diff(S^1)$ by composition from the \emph{right} transforms
coordinates in the domain space $S^1$ as $HQ = \Emb(H\cdot S^1,
\mathbb{R}^2)$, while keeping the curve fixed in the range space
$\mathbb{R}^2$. \medskip 

The present paper discussed how the left action of
$\Diff(\mathbb{R}^2)$ and the right action of $\Diff(S^1)$ on
$Q=\Emb(S^1, \mathbb{R}^2)$ may complement each other in formulating a
variational approach for registration of curves under large
deformations. The left action of $\Diff(\mathbb{R}^2)$ corresponds to
deforming the curve by a time-dependent transformation of the
coordinate system in which it is embedded, while leaving the
parameterisation of the curve invariant. The dynamics of this
deformation is formulated as an Euler-Poincar\'e equation for
$\mathbf{J}_{\rm Sing}\in \mathfrak{X} (\mathbb{R}^2)^{\ast}$ that
results in Hamilton's canonical equations for the momentum and
position variables of the curve that comprise the singular-solution
momentum map (\ref{sing-soln}).  This solution provides the dynamics
for curves that fulfills D'Arcy Thompson's vision of shape
transformation \cite{Th1917}. This vision underlies common practice in
image registration \cite{beg-thesis}.  \medskip

The right action of $\Diff(S^1)$ corresponds to adaptively
reparameterising the $S^1$ domain coordinates of the curve. This
reparameterisation of the curve could be useful, for example, in
designing numerical methods that enhance the resolution of its main
features as it deforms in $\mathbb{R}^2$. As we have seen, this right
action unlocks the parameterisation in the control problem to allow it
more freedom for matching the curve shapes using an equivalent
boundary value problem, without being constrained to match
corresponding points along the template and target curves at the
endpoint in time. As explained above, the action of $\Diff (S^1)$ from
the right gives us the momentum map $\mathbf{J}_S: T ^{\ast} \Emb(S^1,
\mathbb{R}^2) \rightarrow \mathfrak{X}(S)^{\ast}$, which we used to
ensure that the momentum of the curve has no tangential
component. This momentum map for right action is given explicitly as
\[
\mathbf{J}_S = \MM{p}\cdot\pp{}{s}\MM{q}\,.
\]
\medskip

The two momentum maps may be assembled into a single figure as in
\cite{HoMa2004}:

\begin{picture}(150,100)(-105,0)%  
\put(100,75){$T^{\ast} \operatorname{Emb}(S^1,\mathbb{R}^2)$} 
%top label

\put(78,50){$\mathbf{J}_{\rm Sing}$}        
%left label

\put(160,50){$\mathbf{J}_S$}   
%right arrow label

\put(72,15){$\mathfrak{X} (\mathbb{R}^2)^{\ast}$}       
%left bottom label

\put(170,15){$\mathfrak{X}(S^1)^{\ast}$}       
%right bottom label

\put(130,70){\vector(-1, -1){40}}  
% left slanted arrow

\put(135,70){\vector(1,-1){40}}  
% right slanted arrow

\end{picture}
%%%%%%%%%%%%%%%%%%%%%%%%%%%%%%%%%%%%%%%

We use the compatibility of these two momentum maps proven in
\cite{HoMa2004} to divide the curve matching problem into independent
registration and reparameterisation problems, leading to the
reformulation of the curve matching problem as an equivalent geodesic
boundary value problem.

We are currently developing numerical algorithms for the transformed
equations applied to embedded curves and surfaces. As noted in
\cite{Vi2009}, applying a piecewise constant representation to the $q$
and $p$ variables in the untransformed equations results in a set of
ordered points. When this approach is extended to the transformed
equations, a finite volume method is obtained, with the extra terms
taking the form of an advection term in the $q$ equation, and a
continuity term in the $p$ equation, which are very well developed in
the finite volume approach. We will also investigate discontinuous
higher-order polynomial representations of $p$ and $q$ which lead to
discontinous Galerkin methods. Since the error in the transformed
problem can be quantified in terms of the Euclidean distance between
points on the curve for each parameter value, the reparameterised
formulation also makes it much easier to perform Bayesian studies in
which one observes points on a curve with observation error from some
probability distribution, and then one attempts to estimate the
probability distribution for the initial conditions of $p$ (and $\nu$)
for which specified points on the curve match the actual position of
the observed points, after applying the time-1 flow map of the
transformed geodesic equations for $p$ and $q$.  This approach provides a considerably simplified observation operator to which algorithms  such as the Monte Carlo Markov Chain algorithm can be applied.

\subsection*{Acknowledgements} We are grateful to D. C. P. Ellis, F. Gay-Balmaz,  T. S. Ratiu, A. Trouv\'e, F.-X. Vialard and L. Younes for valuable discussions. We thank the Royal Society of London Wolfson Award Scheme for partial support during the course of this work.

\bibliography{reparam}

\begin{thebibliography}{BCHM00}

\bibitem[BCHM00]{BlCrHoMa2000}
A.~M. Bloch, P.~E. Crouch, D.~D. Holm, and J.~E. Marsden.
\newblock An optimal control formulation for inviscid incompressible fluid
  flow.
\newblock In {\em Proc. CDC IEEE}, volume~39, pages 1273--1279, 2000.

\bibitem[BCMR98]{BlCrMaRa1998}
A.~M. Bloch, P.~E. Crouch, J.~E. Marsden, and T.~S. Ratiu.
\newblock Discrete rigid body dynamics and optimal control.
\newblock In {\em Proc. CDC IEEE}, 1998.

\bibitem[Beg03]{beg-thesis}
M.~F. Beg.
\newblock {\em Variational and Computational Methods for Flows of
  Diffeomorphisms in Image Matching and Growth in Computational Anatomy}.
\newblock PhD thesis, Johns Hopkins University, 2003.

\bibitem[CH09]{CoHo2009}
C.~J. Cotter and D.~D. Holm.
\newblock Continuous and discrete {C}lebsch variational principles.
\newblock {\em Found. Comput. Math.}, 9(2):221--242, 2009.

\bibitem[Cot08]{Co2008}
C.~J. Cotter.
\newblock The variational particle-mesh method for matching curves.
\newblock {\em J. Phys. A: Math. Theor.}, 41:344003, 2008.

\bibitem[CY01]{CaYo01}
V.~Camion and L.~Younes.
\newblock Geodesic interpolating splines.
\newblock In M.~Figueiredo, J.~Zerubia, and A.~K. Jain, editors, {\em EMMCVPR
  2001}, volume 2134 of {\em Lecture {N}otes in {C}omputer {S}ciences}.
  Springer, 2001.

\bibitem[GTY04]{GlTrYo04}
J.~Glaunes, A.~Trouv{\'e}, and L.~Younes.
\newblock Diffeomorphic matching of distributions: {A} new approach for
  unlabelled point-sets and sub-manifolds matching.
\newblock In {\em IEEE Computer Society Conference on Computer Vision and
  Pattern Recognition}, volume~2, pages 712--718, 2004.

\bibitem[HM04]{HoMa2004}
D.~D. Holm and J.~E. Marsden.
\newblock Momentum maps and measure valued solutions of the {Euler-Poincar\'e}
  equations for the diffeomorphism group.
\newblock {\em Progr. Math.}, 232:203--235, 2004.
\newblock \url{http://arxiv.org/abs/nlin.CD/0312048}.

\bibitem[Hol09]{Ho2009}
D.~D. Holm.
\newblock {\em Geometric Mechanics, Part II: Translating and Rolling},
  chapter~1.
\newblock Imperial College Press, London, 2009.

\bibitem[MM06]{McMa2006}
R.~I. McLachlan and S.~Marsland.
\newblock The {K}elvin-{H}elmholtz instability of momentum sheets in the
  {E}uler equations for planar diffeomorphisms.
\newblock {\em SIAM Journal on Applied Dynamical Systems}, 5(4), 2006.

\bibitem[MMS06]{MiMaSh2006}
A.~Mills, S.~Marsland, and T.~Shardlow.
\newblock {\em Biomedical Image Registration}, chapter Computing the Geodesic
  Interpolating Spline, pages 169--177.
\newblock Springer, 2006.

\bibitem[She94]{Sh1994}
J.~R. Shewchuk.
\newblock An introduction to the conjugate gradient method without the
  agonizing pain.
\newblock
  \url{www.cs.cmu.edu/\~{}quake-papers/painless-conjugate-gradient.pdf}, 1994.

\bibitem[Tho17]{Th1917}
D'A. Thompson.
\newblock {\em On Growth and Form}.
\newblock 1917.

\bibitem[VG05]{VaGl2005}
Marc Vaillant and Joan Glaunes.
\newblock Surface matching via currents.
\newblock In {\em IPMI}, pages 381--392, 2005.

\bibitem[Via09]{Vi2009}
F.-X. Vialard.
\newblock {\em Approche {H}amiltonienne Pour les Espaces de Formes Dans le
  Cadre des Diff\`eomorphismes: Du Probl\`eme de Recalage d'Images Discontinues
  \`a un Mod\`ele Stochastique de Croissance de Formes}.
\newblock PhD thesis, Ecole Normale Sup\'ereure de Cachan, 2009.

\end{thebibliography}
\end{document}